\documentclass[12pt,reqno]{article}
\usepackage{amsfonts}
\usepackage[utf8x]{inputenc}
\usepackage[english]{babel}
\usepackage{amsmath,amssymb,amsfonts,wasysym}   
\usepackage{graphics}                           
\usepackage{amsmath}  
\usepackage{mathtools} 
\usepackage{hyperref}             
\usepackage{comment}
\newtheorem{theorem}{Theorem}[section]

\newcommand*{\uu}{\mathop{}\!\mathrm{u}}
\newtheorem{proposition}[theorem]{Proposition}

\newtheorem{observation}[theorem]{Observation}
\newtheorem{example}[theorem]{Example}

\newcommand*{\pp}{\mathop{}\!\mathrm{pp}}
\newcommand*{\erf}{\mathop{}\!\mathrm{erf}}

\newcommand*{\vv}{\mathop{}\!\mathrm{v}}

\newcommand*{\sech}{\mathop{}\!\mathrm{sech}}
\newcommand*{\sgn}{\mathop{}\!\mathrm{sgn}}

\newcommand*{\gd}{\mathop{}\!\mathrm{gd}}
\newcommand*{\leb}{\mathop{}\!\mathrm{Leb}}
\newcommand*{\agm}{\mathop{}\!\mathrm{AGM}}
\newcommand*{\tr}{\mathop{}\!\mathrm{Tr}}
\newcommand*{\Proof}{\mathop{}\!\mathit{Proof}}
\newcommand*{\id}{\mathop{}\!\mathrm{id}}

\usepackage{thmtools}
\usepackage{thm-restate}
\newcommand{\textoverline}[1]{$\overline{\mbox{#1}}$}\usepackage{hyperref}

\usepackage{cleveref}

\begin{document}
\title{Gaussian integrals depending by a quantum parameter in finite dimension}
\author{Simone Camosso$^*$}
\date{}
{\renewcommand{\thefootnote}{\fnsymbol{footnote}}
\setcounter{footnote}{1}
\footnotetext{\textbf{e-mail}: r.camosso@alice.it, Saluzzo (CN)}
\setcounter{footnote}{0}
}

\maketitle

\begin{center}
\textbf{Abstract}
\end{center}
\noindent
A common theme in mathematics is the evaluation of Gauss integrals. This, 
coupled with the fact that they are used in different branches of science, makes 
the topic always actual and interesting. 
In these notes we shall analyze a particular class of Gaussian integrals that 
depends by the quantum parameter $\hbar$. Starting from classical results, we will 
present an overview on methods, examples and analogies regarding the practice of 
solving quantum Gaussian integrals.

\vspace{1cm}


\tableofcontents
\section{Introduction}
Let $\hbar$ be a quantum parameter, an important practice in quantum 
mechanics, is the evaluation of Gaussian integrals. 
We will focus our attention on a group of Gaussian integrals that come 
from quantum mechanics and quantum field theory (QFT). A general Gaussian integral in this 
class depends by a quantum parameter $k\,=\,\frac{1}{\hbar}$ with `` $\hbar \rightarrow 0$ ''. 
In our calculations we shall treat $k\,=\,\frac{1}{\hbar}$ as a purely formal 
parameter. Summarising, we will treat integrals of the following form:

\begin{equation}
\label{gaussclassgenint}
\int_{\mathbb{R}^n}A\left(v,w\right)e^{-\frac{i}{\hbar}\psi(v,w)}dv,
\end{equation}
\noindent
where $A$ is the ``amplitude'', $\psi$ is generally quadratic and $v,w 
\in\mathbb{R}^{n}$.

This paper is based on previous discussions of \cite{nove}, \cite{dieci}, \cite{uno} and \cite{sei} on different 
Gaussian integrals of the form $(\ref{gaussclassgenint})$. A basic text is the work of \cite{tre} and, for 
the asymptotic analysis we treat only the Laplace method described in \cite{otto}. 
In order to get the result of a Gaussian integral, sometimes it is necessary to 
introduce a special function (a good reference for the properties of special functions is \cite{ventuno}). Instead on the techniques used, an exhaustive source is the book of Nahin 
\cite{due}. 

Regarding the different examples discussed in these pages, they have been 
inspired by the following works concerning the geometric quantization, they are \cite{quindici}, \cite{sette},
\cite{tredici}, \cite{quattordici} and \cite{undici}.

At the beginning of the paper we will review basic facts concerning Gaussian 
integrals in $1$ and $n$--dimension. 
We proceed with the definition of Gaussian matrix integrals and with different examples coming from geometric quantization.

A special section of this work is dedicated to the connection between Gaussian
integrals and elliptic functions. In this special section old results of number theory 
illustrated by \cite{sedici} and \cite{diciotto}, have been used 
in order to express certain Gaussian integrals in term of the 
arithmetic--geometric mean $\agm$.

The next section is dedicated to the Berezin--Gaussian integral. The Berezin integral is defined for G--numbers 
(Grassmannian numbers) that are noncommutative quantities. The generalized version of the Gaussian 
integral for these variables presents interesting mathematical properties (details are in \cite{nove}).

The last section concerns the relation between Gaussian integrals and Boys integrals. Boys integrals are 
used in quantum chemistry and they are of fundamental importance for many electron systems (\cite{ventidue},\cite{ventitre}).

In the end we want to specify that all examples and discussion has been 
considered only for the finite dimensional case. There are generalizations of 
the Gaussian integral in infinite dimension. For example it is possible consider 
Gaussian integrals in Hilbert spaces, but the discussion of this subject is not part 
of these notes (examples of interesting articles on this suject are \cite{trentaquattro} and \cite{trentacinque}).
The reader can think about this paper as a manual (sure not complete!) or a sort of digest for common use 
by both, students and teachers.

\section{The Gaussian integral in $1$--dimension}

We start with the Gauss integral:

\begin{equation}
\label{gauss1}
\int_{-\infty}^{+\infty}e^{-x^2} dx \, = \,\sqrt{\pi}.
\end{equation}

The result follows taking the square of the left--hand side and using the polar coordinates:

\begin{equation}
\label{Gauss11}
\begin{multlined}[t][12.5cm]
\left(\int_{-\infty}^{+\infty}e^{-x^2} dx\right)^2\,=\, \int_{-\infty}^{+\infty}\int_{-\infty}^{+\infty}e^{-x^2-y^2} dxdy  \\
\,=\,\int_{0}^{2\pi}\int_{0}^{+\infty}\rho e^{-\rho^2} d\rho d\vartheta \,=\, 2\pi\int_{0}^{+\infty}\rho e^{-\rho^2} d\rho\\
\,=\, \pi\int_{0}^{+\infty}e^{-s}ds\,=\, \left. \pi \cdot\lim_{t\rightarrow +\infty}-\frac{1}{e^s}\right|^{t}_{0} \,=\, \pi. 
\end{multlined}
\end{equation}

In his article (\cite{ventiquattro}) Conrad gives eleven proofs of the integral $(\ref{gauss1})$
using different methods.
A quick method consist to use the non elementary formula:

\begin{equation}
 \label{nonelementaryf}
 \Gamma(s)\cdot \Gamma(1-s)\,=\, \frac{\pi}{\sin{(\pi s)}},
\end{equation}
\noindent 
where $\Gamma(s)\,=\, \int_{0}^{+\infty}e^{-x}x^{s-1}dx$ (for $s<0$). Using $(\ref{nonelementaryf})$ 
we have that:

\begin{equation}
\label{alternativeprimenumbermethods}
\begin{multlined}[t][12.5cm]
\left(\int_{-\infty}^{+\infty}e^{-x^2} dx\right)^2\,=\, \left(2\int_{0}^{+\infty}e^{-x^2}dx\right)^2 \\ \,=\, \Gamma\left(\frac{1}{2}\right)\cdot \Gamma\left(\frac{1}{2}\right)\,=\,\frac{\pi}{\sin{\left(\frac{\pi}{2}\right)}}\,=\, \pi. 
\end{multlined}
\end{equation}

The curious reader can read the interesting article of \cite{venticinque} where the author 
treat the integrals as a ``puzzle'' to solve.

Another Gaussian integral is the following:

\begin{equation}
\label{gauss2}
\int_{-\infty}^{+\infty}e^{-ax^2} dx \, = \,\sqrt{\frac{\pi}{a}},
\end{equation}
\noindent
for every $a>0$. We can modify the original integral in order to obtain different versions. For example adding the term $-bx$ we find that:

\begin{equation}
\label{gauss3}
\int_{-\infty}^{+\infty}e^{-ax^2-bx} dx \, = \,e^{\frac{b^2}{4a}}\sqrt{\frac{\pi}{a}}.
\end{equation}

We can prove the formula completing the square $-ax^2-bx\,=\, - \left(\sqrt{a}x +\frac{b}{2\sqrt{a}}\right)^2+\frac{b^2}{4a}$. In this case: 

\begin{equation}
\label{Gauss33}
\begin{multlined}[t][12.5cm]
\int_{-\infty}^{+\infty}e^{-ax^2-bx} dx \,=\, \int_{-\infty}^{+\infty}e^{-\left(\sqrt{a}x +\frac{b}{2\sqrt{a}}\right)^2+\frac{b^2}{4a}} dx\\
\,=\,e^{\frac{b^2}{4a}}\int_{-\infty}^{+\infty}e^{-\left(\sqrt{a}x +\frac{b}{2\sqrt{a}}\right)^2} dx
\,=\, \frac{e^{\frac{b^2}{4a}}}{\sqrt{a}}\int_{-\infty}^{+\infty}e^{-s^2} ds\,=\,e^{\frac{b^2}{4a}}\sqrt{\frac{\pi}{a}}. 
\end{multlined}
\end{equation}

Alternatives of the Gaussian integral are in the book of Nahin \cite{due}, an example is:

\begin{equation}
\label{gauss4}
\int_{0}^{+\infty}x^{2n}e^{-x^2}dx\,=\, \frac{(2n)!\sqrt{\pi}}{2 n! 4^n},
\end{equation}
\noindent
for $n\geq 0$, where for $n\,=\,0$ we refind $(\ref{gauss1})$ between $0$ and $+\infty$. Another version of $(\ref{gauss4})$ is the following:

\begin{equation}\label{gaussvenutodopoahahah} \int_{-\infty}^{+\infty}x^{n}e^{-a x^2}dx\,=\, \frac{1\cdot 3\cdot 5 \cdots 
(n-1)\sqrt{\pi}}{2^{\frac{n}{2}}a^{\frac{n+1}{2}}}, \ \ \ \ n=2,4,\ldots\end{equation}
\noindent
where $a>0$, or this:

\begin{equation}
\label{gauss5}
\int_{0}^{+\infty}\frac{e^{-px^2}-e^{-qx^2}}{x^2}dx\,=\, \sqrt{\pi}\left(\sqrt{q}-\sqrt{p}\right),
\end{equation}
\noindent
for $q>p\geq 0$. 

A complex modification of the Gaussian integral $(\ref{gauss3})$ is the following:

\begin{equation}
\label{gauss6}
\int_{-\infty}^{+\infty}e^{-ax^2-\eta x} dx \, = \,e^{\frac{\eta^2}{4a}}\sqrt{\frac{\pi}{a}},
\end{equation}
\noindent
where $\eta\in\mathbb{C}$ and the integral converges uniformly in any compact region. 
Now, the integral defines an analytic function that may be evaluated by taking  $\eta$ to be real and 
then using analytic continuation. To prove this we use the same trick of the integral $(\ref{gauss3})$ and 
the result is true for all $\eta \in \mathbb{C}$. As observed in $\cite{uno}$ when $\eta\,=\, i\xi$ we get that:

\begin{equation}
\label{fourier}
\int_{-\infty}^{+\infty}e^{-ax^2-i\xi x} dx \, = \, \sqrt{\frac{\pi}{a}}\cdot e^{-\frac{\xi^2}{4a}},
\end{equation}
\noindent
or, in other words, the Fourier transform of $e^{-ax^2}$ is the Gaussian $e^{-\frac{\xi^2}{4a}}$. 

Considering the properties of the Fourier transform we can prove the following result.

\begin{proposition}
Let $m$ be a positive integer and $k$ be a real positive constant, then

\begin{equation}
\int_{-\infty}^{+\infty}x^m e^{-i\xi x -\frac{1}{2}kx^2}dx\,=\, \sqrt{2\pi}\frac{(-i)^m}{k^{m+\frac{1}{2}}}P_{m}(\xi)e^{-\frac{1}{2k}\xi^2},
\end{equation}
\noindent
where $P_{m}(\xi)\,=\,\xi^m+\sum_{j\geq 1}p_{mj}\xi^{m-2j}$ is a monic polynomial in $\xi$ of degree $m$ and parity $(-1)^m$.
\end{proposition}

$\Proof.$\\

By a property of the Fourier transform:

$$ \mathcal{F}\left(x^me^{-\frac{1}{2}kx^2}\right)\,=\,i^m\frac{d^m}{d\xi^m}\mathcal{F}\left(e^{-\frac{1}{2}kx^2}\right).$$

Now the Fourier transform of $e^{-\frac{1}{2}kx^2}$ is equal to $\sqrt{\frac{2\pi}{k}}e^{-\frac{\xi^2}{2k}}$ 
and deriving $m$ times $e^{-\frac{1}{2k}\xi^2}$ we find the polynomial $P_{m}$ where the principal term 
has coefficient $\frac{(-1)^m}{k^m}$. We must collect this coefficient in order to find the polynomial $P_{m}$.

\hfill $\Box$

In one dimension we can do much more using special functions as the gamma function and its inductive property $\Gamma(n+1)\,=\,n\cdot\Gamma(n)$. 
For example we have the cubic Gauss integral: 

\begin{equation}
\label{gauss7}
\int_{0}^{+\infty}e^{-x^3} dx \, = \, \Gamma\left(\frac{4}{3}\right).
\end{equation}

This can be proved by a simple substitution $y\,=\, x^3$:

\begin{equation}
\label{Gauss33}
\begin{multlined}[t][12.5cm]
\int_{0}^{+\infty}e^{-x^3} dx \,=\, \frac{1}{3}\int_{0}^{+\infty}e^{-y} y^{-\frac{2}{3}}dy \,=\,  \frac{1}{3}\cdot\Gamma\left(\frac{1}{3}\right)\,=\, \Gamma\left(\frac{4}{3}\right).
\end{multlined}
\end{equation}

This last argument can be generalized for integrals of the form:

\begin{equation}
\label{gauss7}
\int_{0}^{+\infty}e^{-x^m} dx \, = \, \Gamma\left(\frac{1+m}{m}\right),
\end{equation}
\noindent 
where $m\geq 1$. We observe that for the case $m\,=\,2$ the gamma function is $\Gamma\left(\frac{3}{2}\right)\,=\,\frac{\sqrt{\pi}}{2}$ and we recover $(\ref{gauss1})$ between $0$ and $+\infty$.

Another approach involves the use of the Taylor expansion. In order to explain the method we consider the following integral:

\begin{equation}
\label{gauss8}
\int_{-\infty}^{+\infty} e^{-e^{-x^{2}}} -1 \,dx.
\end{equation}

In this case we can work ``formally'' and write that:

 $$  e^{-e^{-x^{2}}} -1 \,=\, -e^{-x^{2}} +\frac{1}{2!}e^{-2x^2}-\frac{1}{3!}e^{-3x^2} + \cdots, $$
 \noindent 
where we use the Mc. Laurin expansion of $e^{x}\sim 1+x+\frac{x^{2}}{2!}+\cdots$. Now integrating each term we have that:

\begin{equation}
\label{gauss8}
\int_{-\infty}^{+\infty} e^{-e^{-x^{2}}} -1 \,dx\,=\, \sum_{k=1}^{+\infty}\frac{(-1)^{k}}{k!}\int_{-\infty}^{+\infty}e^{-kx^2}dx \,=\,\sum_{k=1}^{+\infty}\frac{(-1)^{k}}{k!}\sqrt{\frac{\pi}{k}}.
\end{equation}
 
So we have:

\begin{equation}
\label{gauss9}
\int_{-\infty}^{+\infty} e^{-e^{-x^{2}}} -1 dx\,=\, \sum_{k=1}^{+\infty}\frac{(-1)^{k}}{k!}\sqrt{\frac{\pi}{k}}.
\end{equation}

The reader can try with $\int_{0}^{+\infty}e^{-e^{-x^{3}}} -1 \, dx$, have a nice exercise!

\subsection{Ways to write the Gaussian integral}

There are different ways to write the Gaussian integral $(\ref{gauss1})$, for example with the substitution $x^2\,=\,\ln{t}$, we have:

\begin{equation}
\label{gaussway1}
\int_{0}^{+\infty} \frac{1}{2\sqrt{\ln{t}}} dt\,=\, \sqrt{\pi}.
\end{equation}

Another substitution is to put $e^x\,=\, s$, so:

\begin{equation}
\label{gaussway12}
\int_{0}^{+\infty} s^{-\ln{s}-1} ds\,=\, \sqrt{\pi}.
\end{equation}

A no obvious way to write the Gaussian integral $(\ref{gauss1})$ is the 
following:

\begin{equation}
\label{gaussway12}
\int_{-\infty}^{+\infty} \frac{1+\sin{\left(\arctan{e^{-x^2}}-\frac{\pi}{2}\right)}}{\cos{\left(\arctan{e^{-x^2}}-\frac{\pi}{2}\right)}}dx\,=\, 
\sqrt{\pi},
\end{equation}
\noindent 
where the Gudermannian function $\gd{(x)}\,=\, \int_{0}^{x}\sech{w}dw$ is used in 
order to prove the integral (see \cite{ventotto} for useful mathematical relations).

Another way to express the Gaussian integral is through the error function:

\begin{equation}
\label{errorfunctionG}
\int_{-\infty}^{+\infty}e^{-t^2}dt\,=\, \sqrt{\pi}\erf{(\infty)},
\end{equation}
\noindent 
where

\begin{equation}
\label{errorfunctionGG}
\erf{(x)}\,=\, \frac{2}{\sqrt{\pi}}\int_{0}^{x}e^{-t^2}dt,
\end{equation}
\noindent 
and it has interesting properties (see \cite{ventinove}). For example we can 
express the `` Plasma dispersion function '' in terms of the $\erf$ function.

The `` Plasma dispersion function '' is defined as:

\begin{equation}
 \label{plasma1}
 D(x)\,=\, \frac{1}{\sqrt{\pi}}\int_{-\infty}^{+\infty}\frac{e^{-t^2}}{t-x}dt,
\end{equation}
\noindent 
for $\Im{(x)}>0$, and its analytic continuation to the rest of the complex $x$ 
plane. Instead of the function $(\ref{plasma1})$ we can consider the `` incomplete 
Plasma dispersion function'':

\begin{equation}
 \label{plasma1}
 D(\nu,x)\,=\, \frac{1}{\sqrt{\pi}}\int_{\nu}^{+\infty}\frac{e^{-t^2}}{t-x}dt,
\end{equation}
\noindent 
for $\Im{(x)}>0$. We obtain $D(x)$ as limit of $D(\nu,x)$ for $\nu\rightarrow -\infty$. The function $D(\nu,x)$ 
satisfy the following differential equation:

 $$ \frac{d}{dx}D(\nu,x)+2xD(\nu,x)\,=\, 
 \frac{1}{\sqrt{\pi}}\frac{e^{-\nu^2}}{\nu-x}-1+\erf(\nu).$$
 
 For references see \cite{trenta}. For $\nu\rightarrow -\infty$, from $\erf(-\infty)\,=\,-1$,  we have that:
 
 $$ D'(x)+2xD(x)\,=\, -2 ,$$
 \noindent 
 that is a linear differential equation with solution:
 
 $$ D(x)\,=\, e^{-x^2}\left[c-2\int_{0}^{x} e^{t^2}dt\right].$$
 
 In order to determine $c$ we observe that in the positive part of the complex 
 plane for $x\rightarrow 0$ we have that:
 
 $$ D(x)\,=\, 
 \frac{1}{\sqrt{\pi}}\pp{\left[\int_{-\infty}^{+\infty}\frac{e^{-t^2}}{t}dt\right]}+i\sqrt{\pi}\,=\, 
 i\sqrt{\pi},$$
 \noindent 
 where $\pp$ is the principal part and the argument of the integral is an odd 
 function, so the integral is zero. 
 
In the end we have that:

$$ D(x)\,=\, e^{-x^2}\left[i\sqrt{\pi}-2\int_{0}^{x}e^{t^2}dt\right]\,=\, 
i\sqrt{\pi}e^{-x^2}\left[1+\frac{2i}{\sqrt{\pi}}\int_{0}^{x}e^{t^2}dt\right]\,=\, 
i\sqrt{\pi}e^{-x^2}\left[1+\erf{(ix)}\right],$$
\noindent 
that is the relation searched.  

The last way to write the Gaussian integral, is given by Ramanujan in \cite{trentadue} using 
the following continued fraction for the integral: 
 
 \begin{equation}
   \label{continuedfrac}
   \int_{0}^{a}e^{-x^2}dx\,=\, \frac{\sqrt{\pi}}{2}-\frac{e^{-a^2}}{2a+\frac{1}{a+\frac{2}{2a+\frac{3}{a+\frac{4}{2a+ \frac{5}{\cdots}}}}}}
 \end{equation}
 
 We give a little sketch of the proof. 
 First we observe that the continued fraction $(\ref{continuedfrac})$ is 
 equivalent to:
 
 \begin{equation}
   \label{continuedfrac2}
   e^{a^2}\int_{a}^{+\infty}e^{-x^2}dx\,=\,\frac{1}{2a+\frac{1}{a+\frac{2}{2a+\frac{3}{a+\frac{4}{2a+ 
   \frac{5}{\cdots}}}}}},
 \end{equation} 
 \noindent
 in fact $\frac{\sqrt{\pi}}{2}\,=\,\int_{0}^{+\infty}e^{-x^2}dx$ so  $\int_{0}^{a}e^{-x^2}dx 
 -\int_{0}^{+\infty}e^{-x^2}dx=-\int_{a}^{+\infty}e^{-x^2}dx$, now we change the sign and 
 multiply by $e^{a^2}$ in order to obtain $(\ref{continuedfrac})$. 
 
Second we consider the generating function:
 
 $$ \psi(x)\,=\, e^{x^2}\int_{x}^{+\infty}e^{-x^2}dx.$$
 
 The function $\psi$ has interesting properties, for example the $n$--th 
 derivative is given by:
 
 $$ \psi^{(n)}(x)\,=\, P_{n}(x)\psi(x)-Q_{n}(x),$$
 \noindent
 where $P_{n}$ and $Q_{n}$ are given by the following recurrence relations:
 
 $$ P_{0}(x)=1, \ \ P_{1}(x)=2x; $$
 $$ Q_{0}(x)=0, \ \ Q_{1}(x)=1; $$
 $$ P_{n+1}(x)=2xP_{n}(x)+2nP_{n-1}(x); $$
 $$ Q_{n+1}(x)=2xQ_{n}(x)+2nQ_{n-1}(x).$$
 
 Furthermore we have the following identity proved by induction:

$$  Q_{n+1}(x)P_{n}(x)-P_{n+1}(x)Q_{n}(x)\,=\,(-2)^nn!$$ 
 \noindent 
 and the following expression for $P_{n}(x)$:
 
 $$ P_{n}(x)\,=\, e^{-x^2}\frac{d^n}{dx^n}e^{x^2}.$$
 
 Third and final step we have to prove that:
 
 $$\lim_{n\rightarrow +\infty}\frac{Q_{n}(x)}{P_{n}(x)}\,=\, \psi(x),$$
 \noindent 
 for all $x>0$. The limit can be showed using upper and lower bounds on $\psi(x)$ 
 in a similar way of \cite{trentatre} that in its article used a different generating 
 function $\phi(x)\,=\,e^{\frac{x^2}{2}}\int_{x}^{+\infty}e^{-\frac{x^2}{2}}dx$ 
 instead our $\psi$. 
 
 Now it is clear that $\frac{Q_{n}(x)}{P_{n}(x)}$ is the $n$--th convergent of 
 the continued fraction:
 
 $$\frac{1}{2a+\frac{2}{2a+\frac{4}{2a+\frac{6}{2a+\frac{8}{2a+ 
   \frac{10}{\cdots}}}}}} $$
   \noindent
that is equivalent simplifying by $2$ to:

     $$\frac{1}{2a+\frac{1}{a+\frac{2}{2a+\frac{3}{a+\frac{4}{2a+ 
   \frac{5}{\cdots}}}}}} $$
   \noindent 
   and the result $(\ref{continuedfrac})$ is proved. Passing to the limit:
   
\begin{equation}
   \label{continuedfraclimit}
   \int_{-\infty}^{+\infty}e^{-x^2}dx\,=\, 2\lim_{a\rightarrow +\infty}\left(\frac{\sqrt{\pi}}{2}-\frac{e^{-a^2}}{2a+\frac{1}{a+\frac{2}{2a+\frac{3}{a+\frac{4}{2a+ 
   \frac{5}{\cdots}}}}}}\right).
 \end{equation}
 
 \subsection{Gaussian integral or Gauss theorem?}
 
The reader that is familiar with physics knows the statement of the Gauss 
theorem for the flux of the electric field $\vec{E}$ generated by a finite number of charges $q_{i}$
confinated in a closed surface. This statement can be summarized by the formula:

\begin{equation}
  \label{electricfieldgauss}
  \Phi\left(\vec{E}\right)\,=\, \frac{\sum_{i=1}^{n}q_{i}}{\varepsilon_{0}},
\end{equation}
\noindent 
where $\varepsilon_{0}$ is called the dielettric constant. 

It is curious that the formula can be expressed in this other complicated way:

\begin{equation}
  \label{electricfieldgauss2}
  \Phi\left(\vec{E}\right)\,=\, \sgn{\left(\sum_{i=1}^{n}q_{i}\right)}\cdot \left(\int_{-\infty}^{+\infty}e^{-\frac{\pi 
  \varepsilon_{0}}{\left|\sum_{i=1}^{n}q_{i}\right|}y^2}dy\right)^2.
\end{equation}

This is a way to express the Gauss theorem with a Gaussian integral! 
This formulation seems not very useful, may be or may be not! 
The reader can imagine heuristically to write the charges $q_{i}$ as an integer 
multiplied by the fundamental charge $e\,=\,1,6\cdot 10^{-19}C$. Now $e$ is 
related to $\hbar$ and it is possible to substitute inside the integral the 
parameter $\hbar$. Somebody can Taylor expand the integral in powers of $\hbar$, 
can this provides a quantization for the flux? 
The interested reader can think about these as formal manipulations, or simply let it go and move on to the next 
section.

\section{The Gaussian integral in $n$--dimensions}

In $n$ dimension we have that:

\begin{equation}
\label{gauss1n}
\int_{\mathbb{R}^{n}}e^{-\|x\|^2} dx \, = \,\sqrt{\pi^{n}},
\end{equation}
\noindent
where $\|x\|^2\,=\, x_{1}^{2}+ \cdots +x_{n}^{2}$. The result follows observing that $\int_{\mathbb{R}^{n}}e^{-\|x\|^2} dx \, = \,\left(\int_{-\infty}^{+\infty} e^{-t^2}dt\right)^{n}$ and using the polar coordinates:

\begin{equation}
\label{Gauss21n}
\begin{multlined}[t][12.5cm]
\frac{1}{\sqrt{\pi^{n}}}\left(\int_{-\infty}^{+\infty}e^{-t^2} dt\right)^n\,=\, \frac{1}{\sqrt{\pi^{n}}}\int_{0}^{+\infty}e^{-\rho^2}\rho^{n-1} c_{n} d\rho \\
\,=\,\frac{c_{n}}{2}\frac{1}{\sqrt{\pi^{n}}}\int_{0}^{+\infty} s^{\frac{n}{2}-1}e^{-s} ds\,=\, \frac{1}{\sqrt{\pi^{n}}}\frac{c_{n}}{2}\Gamma\left(\frac{n}{2}\right)\,=\, 1, 
\end{multlined}
\end{equation}
\noindent 
where $c_{n}\,=\, \frac{2\pi^{\frac{n}{2}}}{\Gamma\left(\frac{n}{2}\right)}$ is the area of the unit sphere. The result has been generalized by \cite{tre} for a symmetric positive definite $n\times n$--matrix $A$:

\begin{equation}
\label{gauss2n}
\int_{\mathbb{R}^{n}}e^{-\langle A x,x\rangle} dx \, = \,\frac{\sqrt{\pi^{n}}}{\sqrt{\det{A}}},
\end{equation}
\noindent
where $\langle \cdot , \cdot \rangle$ is the inner product on $\mathbb{R}^{n}$. 
This result is the analogue of $(\ref{gauss2})$ in one dimension. Always in \cite{tre} we can 
find another result for the Fourier transform. We recall here the theorem.

\begin{theorem}[H\"{o}rmander]
If $A$ is a non--singular symmetric matrix and $\Re{A}\geq 0$ the Fourier transform of 
$u(x)\,=\, e^{-\frac{1}{2}\langle Ax,x\rangle}$ is a Gaussian function $\widehat{u}(\xi)\,=\, (2\pi)^{\frac{n}{2}}(\det{B})^{\frac{1}{2}}e^{-\frac{1}{2}\langle B\xi,\xi\rangle}$ 
where $B=A^{-1}$ and the square root is well defined. If $A\,=\, -iA_{0}$ where $A_{0}$ is real, symmetric and non-singular then 
$\widehat{u}(\xi)\,=\, (2\pi)^{\frac{n}{2}}|\det{A_{0}}|^{-\frac{1}{2}}e^{\frac{\pi i \sgn{A_{0}}}{4}-\frac{1}{2}i\langle A_{0}^{-1}\xi,\xi\rangle}$.
\end{theorem}

In the previous theorem the term $\sgn{A_{0}}$ is called the signature of $A_{0}$.

The Gaussian integral $(\ref{gauss2n})$ admits different generalizations, for 
example we can consider the problem to evaluate the integral:

\begin{equation}
\label{gauss2nwicktheorem}
I\,=\,\int_{\mathbb{R}^{n}}x_{i}x_{j}e^{-\frac{1}{2} x^{T}\cdot Ax} dx,
\end{equation}
\noindent 
with $A$ a real symmetric $n\times n$ matrix, $T$ denotes the transpose and $\cdot$ the ordinary product of matrices. 

A general procedure in order to solve Gauss integrals as $(\ref{gauss2nwicktheorem})$ 
consists to introduce a generating function $\mathcal{Z}(J)$ depending by a parameter 
$J$:

$$ J\,=\,\left(\begin{array}{c} J_{1} \\ \vdots \\ J_{n}\end{array}\right),$$
\noindent 
where

$$ \mathcal{Z}(J)\,=\, \int_{\mathbb{R}^n}e^{-\frac{1}{2}x^{T}\cdot Ax+x^{T}\cdot 
J}.$$

Now we have that:

$$ I\,=\, \left. \frac{\partial^2\mathcal{Z}(J)}{\partial J_{i}\partial 
J_{j}}\right|_{J=0}.$$

We will use this version of the Wick's theorem.

\begin{theorem}
  
  \begin{equation}
   \left.\frac{\partial^n}{\partial J_{i_{1}}\cdots \partial 
J_{i_{n}}} \left(e^{\frac{1}{2}J^{T}\cdot A^{-1}J}\right)\right|_{J=0}\,=\, \sum \left(A^{-1}\right)_{i_{p_{1}}i_{p_{2}}}\cdots\left(A^{-1}\right)_{i_{p_{n-1}}i_{p_{n}}},
\end{equation}
\noindent
where the sum is taken over all pairings: $(i_{p_{1}},i_{p_{2}}), \ldots 
(i_{p_{n-1}},i_{p_{n}})$, of $i_{1}, \ldots, i_{n}$.

\end{theorem}

We will use this theorem in the case where $i_{1}=i$ and $i_{2}=j$. We observe 
that there is another way to write $\mathcal{Z}(J)$:

$$ \mathcal{Z}(J)\,=\, \int_{\mathbb{R}^n}e^{-\frac{1}{2}(x-A^{-1}J)^TA(x-A^{-1}J)+\frac{1}{2}J^TA^{-1}J}.$$

If we substitute $y=x-A^{-1}J$ we have that:

$$\mathcal{Z}(J)\,=\, 
e^{\frac{1}{2}J^TA^{-1}J}\int_{\mathbb{R}^n}e^{-\frac{1}{2}y^TAy}dy\,=\, \frac{(2\pi)^{\frac{n}{2}}}{\sqrt{\det{A}}}e^{\frac{1}{2}J^TA^{-1}J}.$$ 

Reconsidering the integral $I$ we have that:

$$I\,=\, \frac{(2\pi)^{\frac{n}{2}}}{\sqrt{\det{A}}}\left.\frac{\partial^2e^{\frac{1}{2}J^TA^{-1}J}}{\partial J_{i}\partial 
J_{j}}\right|_{J=0}.$$

After calculations:

\begin{equation}
\label{Wick_mega_strange_integral_fuckyou}
\begin{multlined}[t][12.5cm]
I\,=\, \left.\frac{(2\pi)^{\frac{n}{2}}}{\sqrt{\det{A}}}\frac{\partial^2e^{\frac{1}{2}(A^{-1})_{kl}J^kJ^l}}{\partial J_{i}\partial 
J_{j}}\right|_{J=0}\\
\,=\, \left.\frac{(2\pi)^{\frac{n}{2}}}{\sqrt{\det{A}}} \frac{\partial \left[(\frac{1}{2}(A^{-1})_{ai}J^a+\frac{1}{2}(A^{-1})_{ib}J^b)e^{\frac{1}{2}(A^{-1})_{kl}J^kJ^l}\right]}{\partial 
J_{j}}\right|_{J=0}\\
\,=\,\left.\frac{(2\pi)^{\frac{n}{2}}}{\sqrt{\det{A}}} \frac{\partial \left[((A^{-1})_{ia}J^a)e^{\frac{1}{2}(A^{-1})_{kl}J^kJ^l}\right]}{\partial 
J_{j}}\right|_{J=0}\\
\,=\, 
\left.\frac{(2\pi)^{\frac{n}{2}}}{\sqrt{\det{A}}}\left[(A^{-1})_{ij}+(A^{-1})_{ia}J^a\left(\frac{1}{2}(A^{-1})_{bj}J^b+\frac{1}{2}(A^{-1})_{jc}J^c\right)\right]e^{\frac{1}{2}(A^{-1})_{kl}J^kJ^l}\right|_{J=0}\\
\,=\, 
\left.\frac{(2\pi)^{\frac{n}{2}}}{\sqrt{\det{A}}}\left[(A^{-1})_{ij}+(A^{-1})_{ia}J^a\cdot(A^{-1})_{bj}J^b\right]e^{\frac{1}{2}(A^{-1})_{kl}J^kJ^l}\right|_{J=0}\\
\,=\, \frac{(2\pi)^{\frac{n}{2}}}{\sqrt{\det{A}}}(A^{-1})_{ij}.
\end{multlined}
\end{equation}

We may call $\mathcal{Z}_{0}=\frac{(2\pi)^{\frac{n}{2}}}{\sqrt{\det{A}}}$ so the 
result is $I=\mathcal{Z}_{0}\left(A^{-1}\right)_{ij}$. This work can be 
generalized and, a good reference is \cite{ventisette}.

\section{Gaussian integrals and homogeneous functions}

A real function $\varphi:\mathbb{R}^n \rightarrow [0,+\infty[$, at least upper 
semicontinuous is $\alpha$--homogeneous, with $\alpha\,=\,\left(\alpha_{1},\ldots ,\alpha_{n}\right)$ 
and $\alpha_{i}>0$ if

\begin{equation}
 \label{homog1}
 \varphi\left(t^{\alpha_{1}}x_{1},\ldots,t^{\alpha_{n}}x_{n}\right)\,=\,t\varphi\left(x_{1},\ldots,x_{n}\right)
\end{equation}
\noindent
for each $x\in\mathbb{R}^n$ and $t>0$.

We may consider the sphere for the $\alpha$--homogeneous function $\varphi$ defined for each $r>0$:

$$ B_{\varphi}(r)\,=\, \left\{x\in\mathbb{R}^n : \varphi(x)<r\right\}.$$

\begin{observation}
 $B_{\varphi}(1)$ is the unit ball and $\leb{B_{\varphi}(r)}$ denotes the 
 Lebesgue measure of the ball $B_{\varphi}(r)$.
\end{observation}
 
In \cite{trentuno} the following theorem has been proved.

\begin{theorem}
 Let $\varphi:\mathbb{R}^n \rightarrow [0,+\infty[$ be an $\alpha$--homogeneous 
 function, then 
 
 \begin{equation}
 \label{homog2}
\int_{\mathbb{R}^n}e^{-\varphi(x)}dx\,=\, \leb{(B_{\varphi})}p!
\end{equation}
 \noindent 
 where $p\,=\, \alpha_{1}+\cdots +\alpha_{n}$ is the weight of $\varphi$. 
 \end{theorem}

The next proposition is an application of a study on homogeneous functions and 
Gaussian integrals. 

\begin{proposition}
Let $c>0$, $p>0$, and let $A$ be an $n\times n$ positive definite symmetric 
matrix, then 

 \begin{equation}
 \label{homog3}
\int_{\mathbb{R}^n}e^{-(c\langle Ax,x\rangle)^p}dx\,=\, 
\left(\sqrt{\frac{\pi}{c}}\right)^n\frac{1}{\sqrt{\det{A}}}\frac{\Gamma\left(\frac{n}{2p}+1\right)}{\Gamma\left(\frac{n}{2}+1\right)}.
\end{equation}
\end{proposition}

\begin{observation}
From the previous proposition we refind that:

$$\int_{0}^{+\infty}e^{-x^3}dx\,=\, 
\frac{1}{2}\int_{-\infty}^{+\infty}e^{-|x|^3}dx\,=\, 
\frac{1}{2}\int_{-\infty}^{+\infty}e^{-(x^2)^{\frac{3}{2}}}dx\,=\, 
\sqrt{\pi}\frac{\Gamma\left(\frac{1}{3}+1\right)}{\Gamma\left(\frac{1}{2}+1\right)}\,=\,
\Gamma\left(\frac{4}{3}\right).$$

\end{observation}

\section{Gaussian matrix integrals}

Gaussian integrals may be also defined on the space of Hermitian matrices. 

Let $\mathcal{H}_{N}$ be the space of all Hermitian matrices $H\,=\,\left\{h_{ij}\right\}_{i=1,\ldots,N,j=1,\ldots,N}$ 
where $h_{ij}\in\mathbb{C}$ and $h_{ij}=\overline{h_{ji}}$. Every matrix $H\in\mathcal{H}_{N}$ 
may be represented as:

$$ H\,=\, \left[\begin{array}{cccc} x_{11} & (x_{12}+iy_{12}) & \cdots & (x_{1N}+iy_{1N}) \\  \overline{(x_{12}+iy_{12})} & x_{22} & \cdots & \vdots \\ \vdots & \vdots & \ddots & \vdots \\  \overline{(x_{1N}+iy_{1N})} & \cdots & \cdots & x_{NN} \end{array}\right]$$
\noindent 
where $x_{ij}, y_{ij} \in \mathbb{R}$ for $1\leq i<j\leq N$ and $x_{ii}\in \mathbb{R}$ 
for $i=1, \ldots, N$. We have that $\mathcal{H}_{N} \sim \mathbb{R}^{N^2}$ and 
we denote by:

$$ d\textbf{Leb}(H)\,=\, \prod_{i=1}^{N} dx_{ii}\prod_{i<j}dx_{ij}dy_{ij},$$
\noindent 
the Lebesgue measure in $\mathcal{H}_{N}$. 

In order to introduce a Gaussian measure we must define a non degenerate 
quadratic form on $\mathcal{H}_{N}$. Let us consider:

$$ H^2\,=\, \sum_{j=1}^{N}h_{ij}h_{jk},$$
\noindent 
and its trace:

\begin{equation}
\label{matrixgaussint}
\begin{multlined}[t][12.5cm]
\tr(H^2)\,=\, \sum_{i,j=1, \ldots ,N}h_{ij}h_{ji}\,=\, \sum_{i,j=1, \ldots, N}h_{ij}\overline{h_{ij}} 
\\
\,=\, \sum_{i=j}x_{ii}^2 +\sum_{i\not=j} x_{ij}^{2}+y_{ij}^{2} \\
\,=\, \sum_{i=1}^{N}x_{ii}^{2}+2\sum_{i<j}\left(x_{ij}^{2}+y_{ij}^{2}\right) 
\,=\, \left(Bh,h\right).
\end{multlined}
\end{equation}

where:

$$ h\,=\, \left[\begin{array}{c}x_{11} \\ \vdots \\ x_{NN} \\ x_{12} \\ \vdots \\ x_{N-1,N} \\ y_{12} \\ \vdots \\ y_{N-1,N}\end{array}\right] 
$$
\noindent 
and 

$$B\,=\, \left[\begin{array}{cccccc} 1&&&&& \\ & \ddots & &&& \\ && 1 &&& \\ &&& 2 && \\ &&&&\ddots & \\ &&&&& 2 \end{array}\right]$$
\noindent 
is an $N^2\times N^2$ matrix, with $\det{B}\,=\, 2^{N^2-N}$. We may define the 
Gaussian measure:

\begin{equation}
\label{gaussianmatrixmeasure}
d\mu(H)\,=\, \frac{1}{(\sqrt{2\pi})^{N^2}}\cdot 
2^{\frac{N^2-N}{2}}e^{-\frac{1}{2}\tr(H^2)}d\textbf{Leb}(H),
\end{equation}
\noindent 
and the Gaussian matrix integral:

\begin{equation}
\label{gaussianmatrixintegralfinalversion}
\int_{\mathcal{H}_{N}}e^{-\frac{1}{2}\tr{(H^2)}}d\textbf{Leb}(H)\,=\, 
\frac{(\sqrt{2\pi})^{N^2}}{2^{\frac{N^2-N}{2}}}.
\end{equation}

A complete and depth discussion on the topic can be found in \cite{ventisei}.

\section{The symplectic structure}

In this section we consider $\mathbb{R}^d$ with $d\,=\, 2n$. We assume on $\mathbb{R}^d$ the presence 
of $J: \mathbb{R}^d\rightarrow \mathbb{R}^d$ such that $J^2\,=\, -\id$, where $\id$ is the identity 
map on $\mathbb{R}^d$ and $J(v)\,=\, iv$ for all $v\in \mathbb{R}^d$. We say that $(\mathbb{R}^d,J)$ is a complex vector 
space also denoted by $\mathbb{R}^d_{J}$. Now if we assume to have $H:\mathbb{R}^d_{J}\times \mathbb{R}^d_{J}\rightarrow \mathbb{C}$ be 
a positive definite hermitian product then $g\,=\, \Re{(H)}$ and $\omega\,=\, -\Im{(H)}$ define 
respectively  a real symmetric scalar product and a symplectic structure on $\mathbb{R}^d_{J}$. 
The reverse is also true. In fact given $g: \mathbb{R}^d_{J}\times \mathbb{R}^d_{J}\rightarrow \mathbb{R}$ 
a real $J$--invariant scalar product if we define $\omega(\cdot , \cdot)\,=\, g(J\cdot,\cdot)$ and 
$H\,=\, g-i\omega$ then $H$ is a positive definite hermitian product. 

Furthermore if we define $g(v,w)\,=\, v^{T}\cdot w$ 
where the exponent $T$ is the transpose and $v,w\in \mathbb{R}^d_{J}$, we have that $g$ is a real 
$J$--invariant scalar product and the structure $\omega$ can be described in term of $g$. 
Recalling that $J^{T}\,=\,-J$ we have:

\begin{equation}
\label{symplectic2}
\omega(v , w)\,=\, g(J v,w), \ \ \omega(v , Jw)\,=\, g(v,w), \ \ \omega(Jv,Jw)\,=\, g(Jv,w),
\end{equation}
\noindent 
and also $\omega$ is $J$--invariant.

\section{The quantum postulates}

In his work $\cite{quattro}$, P. Dirac defines the quantum Poisson bracket $[\cdot,\cdot]$ of any two variables $\uu$ and $\vv$ as

\begin{equation}
\label{dirac}
\uu\vv-\vv\uu=i\hbar [\uu,\vv],
\end{equation}
\noindent
where $\hbar\,=\, 6,626070040(81)\cdot 10^{-34} J\cdot s$ is the Planck constant. The formula $(\ref{dirac})$ is one of the basic postulates of quantum mechanics. 
We can summarize these postulates as follows. To start we fix a symplectic manifold $(M,\omega)$ of dimension $d_{M}$, 
with $\omega$ the corresponding symplectic structure and an Hilbert space $\mathcal{H}$. The quantization is a ``way'' to pass from the classical system to the quantum system. 
In this case the classical system, or phase space, is described by the symplectic manifold $M$ and the Poisson algebra of smooth functions on $M$ denoted by $(\mathcal{C}^{\infty}(M),\{\cdot,\cdot\})$. 
The quantum system is described by $\mathcal{H}$. We define quantization a map $Q$ from the subset of 
the commutative algebra of observables $\mathcal{C}^{\infty}(M)$ to the space of operators in $\mathcal{H}$. 
Let $f\in \mathcal{C}^{\infty}(M)$ we have that $Q(f):\mathcal{H}\rightarrow \mathcal{H}$ is the corresponding quantum operator.
We can summarize the quantum axioms in this scheme:

\begin{itemize}
 \item[1.] linearity: $Q(\alpha f+\beta g)\,=\,\alpha Q(f)+\beta Q(g)$, for every $\alpha,\beta$ scalar and $f,g$ observables;
 \item[2.] normality: $Q(1)\,=\,I$, where $I$ is the identity operator;
 \item[3.] hermiticity: $Q(f)^{*}\,=\,Q(f)$;
 \item[4.] (Dirac) quantum condition: $[Q(f),Q(g)]\,=\,-i\hbar Q(\{f,g\})$;
 \item[5.] Irreducibility condition: for a given set of observables $\{f_{j}\}_{j\in I}$, with the property that for every other $g\in\mathcal{C}^{\infty}(M)$, 
 such that $\{f_{j},g\}\,=\,0$ for all $j$, then $g$ is constant. We can associate a set of quantum operators $\{Q(f_{j})\}_{j\in I}$ such that for every other operator
 $Q$ that commute with all of them is a multiple of the identity.
\end{itemize}

The last postulate says that in the case we consider a connected Lie group $G$ we say that it is a group 
of symmetries of the physical system if we have the two following irreducible representation: 
one as symplectomorphisms acting on $(M,\omega)$ and another as unitary transformations acting on $\mathcal{H}$. 
For many details about these postulates see $\cite{cinque}$.

\begin{example}[Schr\"{o}dinger quantization]
Let $M\,=\,\mathbb{R}^{2n}$ and $(q_{j},p_{j})$ the canonical coordinates of position and momentum. 
In this case $Q(q_{j})\,=\,q_{j}$ that acts as multiplication and $Q(p_{j})\,=\,-i \hbar\frac{\partial}{\partial q_{j}}$. 
The space $\mathcal{H}=L^{2}(\mathbb{R}^{n},dq_{j})$ is the quantum space and there are the following relations of commutations:

$$ [Q(q_{k}),Q(q_{j})]=[Q(p_{k}),Q(p_{j})]\,=\,0, \ \ \ \ [Q(q_{k}),Q(p_{j})]\,=\,i\hbar\delta_{kj}I.$$
\end{example}

Quantization models in the past years have been studied from different point of view. We cite as examples
the case of geometric quantization and deformation quantization but the literature is very wide on the topic. 
We will consider the basic case on $\mathbb{R}^{2n}$ presented in the above example 
as quantization scheme for the next sections.

\section{Gaussian integrals depending by a quantum parameter}

In this section we study a particular kind of Gaussian integrals deriving from quantum mechanics. 
The general form of these integrals is the following:

\begin{equation}
 \label{quantgauss}
 \int_{\mathbb{R}^{n}}e^{-\frac{i}{\hbar}g(v,w)}dv,
\end{equation}
\noindent 
where the function $g$ is usually quadratic in its variables and $(v,w)\in \mathbb{R}^{2n}$. 

Here we are interested to particular forms of $g$.  

\begin{proposition}
We have the following cases:

\begin{itemize}
\item[$1)$] if $g\,=\,g_{1}(v,w)\,=\,\omega(v,w)-\frac{i}{2}\|v-w\|^2$ then:

\begin{equation}
\label{result1}
\int_{\mathbb{R}^{n}}e^{-\frac{i}{\hbar}g_{1}(v,w)}dv\,=\,\left(2\pi\hbar\right)^{\frac{n}{2}}e^{-\frac{1}{2\hbar}\|w\|^2};
\end{equation}
\item[$2)$] if $g\,=\,g_{2}(v,w,u)\,=\,-\omega(v,w+u)-\frac{i}{2}\|v\|^2$ then:

\begin{equation}
\label{result2}
\int_{\mathbb{R}^{n}}e^{-\frac{i}{\hbar}g_{2}(v,w,u)}dv\,=\,(2\pi\hbar)^{\frac{n}{2}}e^{-\frac{1}{2\hbar}\|w+u\|^2};
\end{equation}
\item[$3)$] if $g\,=\,g_{3}(v,w,u)\,=\,\omega(v,w)+\omega(w,u)-\frac{i}{2}\|v-w\|^2-\frac{i}{2}\|w-u\|^2$ then:

\begin{equation}
\label{result3}
\int_{\mathbb{R}^{n}}e^{-\frac{i}{\hbar}g_{3}(v,w,u)}dw\,=\,(\pi\hbar)^{\frac{n}{2}}e^{-\frac{i}{\hbar}g_{1}(v,u)},
\end{equation}

\end{itemize}
\noindent
for all $v,w,u\in\mathbb{R}^{n}$.
\end{proposition}

$\Proof$.\\

We start with the first Gaussian integral and we have that:

\begin{equation}
\label{Gaussquant1}
\begin{multlined}[t][12.5cm]
\int_{\mathbb{R}^{n}}e^{-\frac{i}{\hbar}g_{1}(v,w)}dv\,=\,\int_{\mathbb{R}^{n}}e^{-\frac{i}{\hbar}\omega(v,w)-\frac{1}{2\hbar}\|v-w\|^2}dv  \\
\,=\,\hbar^{\frac{n}{2}}\int_{\mathbb{R}^{n}}e^{-ig\left(J(\beta),\frac{w}{\sqrt{\hbar}}\right)-\frac{1}{2}\|\beta\|^2}d\beta\,=\,\left(2\pi\hbar\right)^{\frac{n}{2}}e^{-\frac{1}{2\hbar}\|w\|^2}, 
\end{multlined}
\end{equation}
\noindent 
where $\beta\,=\,\frac{1}{\sqrt{\hbar}}(v-w)$ is a new variable used in the integration. 

For the second case we proceed in a similar way:

\begin{equation}
\label{Gaussquant12}
\begin{multlined}[t][12.5cm]
\int_{\mathbb{R}^{n}}e^{-\frac{i}{\hbar}g_{2}(v,w,u)}dv\,=\,\int_{\mathbb{R}^{n}}e^{\frac{i}{\hbar}\omega(v,w+u)-\frac{1}{2\hbar}\|v\|^2}dv  \\
\,=\,\int_{\mathbb{R}^{n}}e^{-\frac{i}{\hbar}g(Jv,w+u)-\frac{1}{2}\left\|\frac{v}{\sqrt{\hbar}}\right\|^2}dv \,=\,\int_{\mathbb{R}^{n}}e^{-i g\left(J\frac{v}{\sqrt{\hbar}},\frac{w+u}{\sqrt{\hbar}}\right)-\frac{1}{2}\left\|\frac{v}{\sqrt{\hbar}}\right\|^2}dv\\
\,=\,\hbar^{\frac{n}{2}}\int_{\mathbb{R}^{n}}e^{-i g\left(J\beta,\frac{w+u}{\sqrt{\hbar}}\right)-\frac{1}{2}\left\|\beta\right\|^2}d\beta\,=\,(2\pi\hbar)^{\frac{n}{2}}e^{-\frac{1}{2\hbar}\|w+u\|^2},
\end{multlined}
\end{equation}
\noindent 
where we used $\beta\,=\,\frac{1}{\sqrt{\hbar}}v$ as a new variable.

In the last case we consider the function $g_{3}(v,w,u)\,=\, \omega\left(v,w\right) -\frac{i}{2}\|v-w\|^2+\omega\left(w,u\right) -\frac{i}{2}\|w-u\|^2$ where $v,w,u\in\mathbb{R}^{n}$. 

\begin{equation}
\label{Gaussquant13}
\begin{multlined}[t][12.5cm]
\int_{\mathbb{R}^{n}}e^{-\frac{i}{\hbar}g_{3}(v,w,u)}dw\,=\,\int_{\mathbb{R}^{n}}e^{-\frac{i}{\hbar}\left[\omega(v,w)+\omega(w,u)\right]-\frac{1}{2\hbar}\left[\|v-w\|^2+\|w-u\|^2\right]}dw. 
\end{multlined}
\end{equation}

We set $v-w\,=\, t$ so:

\begin{equation}
\label{Gaussquant13newv}
\begin{multlined}[t][12.5cm]
\int_{\mathbb{R}^{n}}e^{-\frac{i}{\hbar}g_{3}(v,w,u)}dw\,=\,e^{-\frac{i}{\hbar}\omega(v,u)}\int_{\mathbb{R}^{n}}e^{-\frac{i}{\hbar}\left[\omega(t,v-u)\right]-\frac{1}{2\hbar}\left[\|t\|^2+\|v-u-t\|^2\right]}dt. 
\end{multlined}
\end{equation}

A second variable $\frac{v-u}{2}-t\,=\,z$ permits to write:

\begin{equation}
\label{Gaussquant13newvv}
\begin{multlined}[t][12.5cm]
\int_{\mathbb{R}^{n}}e^{-\frac{i}{\hbar}g_{3}(v,w,u)}dw\,=\,e^{-\frac{i}{\hbar}\omega(v,u)}\int_{\mathbb{R}^{n}}e^{-\frac{i}{\hbar}\left[\omega(z,u-v)\right]-\frac{1}{2\hbar}\left[\left\|\frac{v-u}{2}-z\right\|^2+\left\|\frac{v-u}{2}+z\right\|^2\right]}dz. 
\end{multlined}
\end{equation}

After calculations we find that:

\begin{equation}
\label{Gaussquant13newvvv}
\begin{multlined}[t][12.5cm]
\int_{\mathbb{R}^{n}}e^{-\frac{i}{\hbar}g_{3}(v,w,u)}dw\,=\,e^{-\frac{i}{\hbar}\omega(v,u)-\frac{1}{4\hbar}\|v-u\|^2}\int_{\mathbb{R}^{n}}e^{-\frac{i}{\hbar}\left[\omega(z,u-v)\right]-\frac{1}{\hbar}\|z\|^2}dz. 
\end{multlined}
\end{equation}

Now the last integral is an ordinary Gaussian integral of simple estimation:

\begin{equation}
\label{Gaussquant13newvvv}
\begin{multlined}[t][12.5cm]
\int_{\mathbb{R}^{n}}e^{-\frac{i}{\hbar}g_{3}(v,w,u)}dw\,=\,(\pi\hbar)^{\frac{n}{2}}e^{-\frac{i}{\hbar}\omega(v,u)-\frac{1}{2\hbar}\|v-u\|^2}\,=\, (\pi\hbar)^{\frac{n}{2}}e^{-\frac{i}{\hbar}g_{1}(v,u)}. 
\end{multlined}
\end{equation}

\hfill $\Box$

A variation of the exponent is the function $g_{4}(v,w)\,=\, \omega(v,Jw)-\frac{i}{2}\|v\|^2-\frac{i}{2}(1+2i\cot{\vartheta})\|w\|^2$. This function appears in \cite{sette} as the function $\psi_{2}(v_{0},e^{i\vartheta}v_{1})$.

\begin{proposition}
Let $\vartheta \in [0,\pi]$, $v,w \in\mathbb{R}^{n}$, then

\begin{equation}
\label{result4}
\int_{\mathbb{R}^{2n}}e^{-\frac{i}{\hbar}g_{4}(v,w)}dvdw\,=\,2^{\frac{n}{2}}(\hbar\pi)^{n}\sin^{\frac{n}{2}}{\vartheta}e^{i\left(\frac{\pi}{2}-\vartheta\right)\frac{n}{2}}.
\end{equation} 
\end{proposition}

$\Proof$.\\
In this case we have that: 

\begin{equation}
\label{Gaussquant2modified}
\begin{multlined}[t][12.5cm]
\int_{\mathbb{R}^{2n}}e^{-\frac{i}{\hbar}g_{4}(v,w)}dvdw\,=\,\int_{\mathbb{R}^{2n}}e^{-\frac{i}{\hbar}\omega(v,Jw)-\frac{1}{2\hbar}\|v\|^2-\frac{1}{2\hbar}(1+2i\cot{\vartheta})\|w\|^2}dvdw  \\
\,=\,\left(\hbar 2\pi\right)^{\frac{n}{2}}\int_{\mathbb{R}^{n}}e^{-\frac{1}{2\hbar}(2+2i\cot{\vartheta})\|w\|^2}dw\,=\, \frac{\hbar^{n}(2\pi)^{n}}{\sqrt{\det{(2+2i\cot{\vartheta})I_{n}}}}, 
\end{multlined}
\end{equation}
\noindent 
where $I_{n}$ is the identity matrix. Observing that $2+2i\cot{\vartheta}\,=\, \frac{2i}{\sin{\vartheta}}\left[-i\sin{\vartheta}+\cos{\vartheta}\right]$ we have that:

\begin{equation}
\label{Gaussquant2modified}
\begin{multlined}[t][12.5cm]
\int_{\mathbb{R}^{2n}}e^{-\frac{i}{\hbar}g_{4}(v,w)}dvdw\,=\, \frac{\hbar^{n}(2\pi)^{n}}{\sqrt{\det{\frac{2ie^{-i\vartheta}}{\sin{(\vartheta)}}I_{n}}}}\,=\,(\hbar\pi)^{n}2^{\frac{n}{2}}\sin^{\frac{n}{2}}{(\vartheta)}e^{i\left(\frac{\pi}{2}-\vartheta\right)\frac{n}{2}}.
\end{multlined}
\end{equation}

\hfill $\Box$

\begin{proposition}
Let $g_{5}(u,v,w)\,=\, \frac{1}{2}(\|w-v\|^{2}+\|v-u\|^{2})$ and  $v,w,u \in\mathbb{R}^{n}$, then

\begin{equation}
\label{result4}
\int_{\mathbb{R}^{n}}e^{\frac{i}{\hbar}g_{5}(u,v,w)}dv\,=\,\left(2\pi\hbar\right)^{\frac{n}{2}}e^{\frac{i}{2\hbar}\left\|\frac{w-u}{2}\right\|^2}e^{\frac{i\pi n}{4}}.
\end{equation} 
\end{proposition}

$\Proof$.\\

From the integral:

\begin{equation}
\label{GaussquantFey}
\begin{multlined}[t][12.5cm]
\int_{\mathbb{R}^{n}}e^{\frac{i}{2\hbar}\left(\|w-v\|^2+\|v-u\|^2\right)}dv\,=\,\int_{\mathbb{R}^{n}}e^{\frac{i}{2\hbar}\left(\|(w-u)-t\|^2+\|t\|^2\right)}dt  \\
\end{multlined}
\end{equation}
\noindent 
where we used the substitution $v-u\,=\, t$. A second substitution $t\,=\,\frac{(w-u)}{2}-z$ bring to the integral:

\begin{equation}
\label{GaussquantFey2}
\begin{multlined}[t][12.5cm]
\,=\,\int_{\mathbb{R}^{n}}e^{\frac{i}{2\hbar}\left(\left\|\frac{(w-u)}{2}+z\right\|^2+\left\|\frac{(w-u)}{2}-z\right\|^2\right)}dz\,=\, \int_{\mathbb{R}^{n}}e^{\frac{i}{2\hbar}\left(\left\|\frac{(w-u)}{2}\right\|^2+\|z\|^2\right)}dz  \\
\,=\,e^{\frac{i}{2\hbar}\left\|\frac{(w-u)}{2}\right\|^2}\int_{\mathbb{R}^{n}}e^{\frac{i}{2\hbar}\|z\|^2}dz\,=\, \hbar^{\frac{n}{2}}e^{\frac{i}{2\hbar}\left\|\frac{(w-u)}{2}\right\|^2}\int_{\mathbb{R}^{n}}e^{\frac{i}{2}\|q\|^2}dq\\
\,=\, \left(2\pi\hbar\right)^{\frac{n}{2}}e^{\frac{i}{2\hbar}\left\|\frac{(w-u)}{2}\right\|^2}e^{\frac{i\pi n}{4}},
\end{multlined}
\end{equation}
\noindent 
where the last substitution was $\frac{1}{\sqrt{\hbar}}z\,=\,q$. 
\hfill $\Box$

\section{Quantum integrals}

Now we shall consider a slight variation in the original integral adding a multiplicative factor $A$ in front of the exponential:

\begin{equation}
 \label{quantgaussgg2}
 \int_{\mathbb{R}^{n}}A(v,w)e^{-\frac{i}{\hbar}\psi(v,w)}dv.
\end{equation}

What is the general approach to this kind of integrals? Without specifying the form 
of the functions $A,\psi$ it is difficult to find a general strategy that works always. 
We can observe that the integral $(\ref{quantgaussgg2})$ is an oscillating integral with parameter $\frac{1}{\hbar}$. 
So the classical approach is the method of stationary phase. 
The idea consists to expand the function $\psi(v,w)$ near its critical points in order to have a ``good'' approximation of the integral. 
References are \cite{tre} and \cite{sei}.

Another very flexible method is the Laplace method. 
It consists to Taylor expand $A,\psi$ around the critical points of $\psi$ and estimate the integral. Let $A,\psi$ be infinitely differentiable functions in $\mathbb{R}^{n}$.
 In our case assume to consider a specific form of $\psi(v,w)=-\frac{i}{2}\|v-w\|^2$. 
 Thus we must evaluate the integral:

\begin{equation}
\label{quantgaussgg3}
I \,=\, \int_{\mathbb{R}^{n}}A(v,w)e^{-\frac{1}{2\hbar}\|v-w\|^2}dv.
\end{equation}

Assuming that $I$ converges absolutely for $\frac{1}{\hbar}$ sufficiently large, we can start changing variable $v-w=u$:

\begin{equation}
\label{quantgaussgg3gg4}
I\,=\,\int_{\mathbb{R}^{n}}A(u+w,w)e^{-\frac{1}{2\hbar}\|u\|^2}du.
\end{equation}

Now the function $\psi(u)\,=\,\|u\|^2$ has $u=0$ as critical point. Expanding $A(u+w,w)$ around zero we have that:

\begin{equation}
\label{quantgaussgg3gg4}
I\,=\,\int_{\mathbb{R}^{n}}A(u+w,w)e^{-\frac{1}{2\hbar}\|u\|^2}du\,=\,\int_{\mathbb{R}^{n}}\left(A(w,w)+DA(u+w,w)|_{u=0}\|u\|+\ldots \right)e^{-\frac{1}{2\hbar}\|u\|^2}du,
\end{equation}
\noindent 
that integrated term by term gives:

\begin{equation}
\label{quantgaussgg3gg4}
\begin{multlined}[t][12.5cm]
I\,=\,\int_{\mathbb{R}^{n}}A(u+w,w)e^{-\frac{1}{2\hbar}\|u\|^2}du\,=\, A(w,w)\left(2\pi \hbar\right)^{\frac{n}{2}}+ \\
+DA(u+w,w)|_{u=0}\frac{2^{\frac{n+1}{2}}\pi^{\frac{n}{2}}\Gamma\left(\frac{n+1}{2}\right)}{\Gamma\left(\frac{n}{2}\right)}\hbar^{\frac{n+1}{2}}+\\+D^2A(u+w,w)|_{u=0}\frac{2^{\frac{n+2}{2}}\pi^{\frac{n}{2}}\Gamma\left(\frac{n+2}{2}\right)}{\Gamma\left(\frac{n}{2}\right)}\hbar^{\frac{n+2}{2}}+\mathcal{O}\left(\hbar^{\frac{n}{2}+1}\right).
\end{multlined}
\end{equation}

We used the Laplace approximation formula in \cite{otto} for $J(k)$ as $k\rightarrow +\infty$:

\begin{equation}
\label{Laplaceapp}
J(k)\,=\,\int_{\mathbb{R}^{n}}A(x)e^{-k\psi(x)}dx\,=\, A\left(\hat{x}\right)e^{-k\psi(\hat{x})}(2\pi)^{\frac{n}{2}}|\Sigma|^{\frac{1}{2}}k^{-\frac{n}{2}}+ \mathcal{O}(1/k),
\end{equation}
\noindent 
where $\Sigma=D^2\psi(\hat{x})^{-1}$ is the Hessian in the critical point $\hat{x}$ of $\psi$. 

Weaker assumptions can be made on the functions $A,\psi$ requiring to be $A$ continuous and $\psi$ has continuous second order derivarites in a neighborhood of $\hat{x}$.

\section{On a Gaussian integral in one dimension and the Lemniscate}

\begin{proposition}
Let $\hbar$ be the Plank constant then:

\begin{equation}
\label{lemniscate1}
\int_{-\infty}^{+\infty}e^{-\frac{1}{4\hbar}x^4+\frac{1}{2\hbar}x^2}dx\,=\, 
\mathcal{C}_{1}(n,\hbar)\cdot \frac{\sqrt{\hbar^{\frac{1}{2}}2^{\frac{1}{2}}\pi^{\frac{3}{2}}}}{\agm(\sqrt{2},1)^{\frac{1}{2}}}.
\end{equation}
\noindent
where

\begin{equation}
\begin{multlined}[t][12.5cm]
\mathcal{C}_{1}(n,\hbar)\,=\,\sum_{n=0}^{+\infty}\frac{1}{\hbar^{\frac{n}{2}}n!}\sqrt{\frac{(2n)!}{2^{3n}n!}
\prod_{k=1}^{+\infty}\frac{\left(k-\frac{1}{2}\right)(2n+4k-1)}{k(2n+4k-3)}
\prod_{l=1}^{+\infty}\frac{2l(4l-3)}{(2l-1)(4l-1)}}.
\end{multlined}
\end{equation}

\end{proposition}

$\Proof.$

\begin{equation}
\begin{multlined}[t][12.5cm]
\int_{-\infty}^{+\infty}e^{-\frac{1}{4\hbar}x^4+\frac{1}{2\hbar}x^2}dx\,=\, \int_{-\infty}^{+\infty}
\left(\sum_{n=0}^{+\infty}\left(\frac{1}{2\hbar}\right)^n\frac{x^{2n}}{n!}\right)e^{-\frac{1}{4\hbar}x^4}dx\\
\,=\,2\int_{0}^{+\infty}
\left(\sum_{n=0}^{+\infty}\left(\frac{1}{2\hbar}\right)^n\frac{x^{2n}}{n!}\right)e^{-\frac{1}{4\hbar}x^4}dx.
\end{multlined}
\end{equation}

We substitute the variable $\frac{x^4}{4\hbar}\,=\,u$, then we have: 

\begin{equation}
\begin{multlined}[t][12.5cm]
\int_{-\infty}^{+\infty}e^{-\frac{1}{4\hbar}x^4+\frac{1}{2\hbar}x^2}dx\,=\, 
\sum_{n=0}^{+\infty}2^{-\frac{1}{2}}\hbar^{\frac{1}{4}}
\int_{0}^{+\infty}\frac{1}{\hbar^{\frac{n}{2}}}\frac{u^{\frac{n}{2}-\frac{3}{4}}}{n!}e^{-u}du,
\end{multlined}
\end{equation}
\noindent
that we can write as $\frac{\hbar^{\frac{1}{4}}}{\sqrt{2}}\sum_{n=0}^{+\infty}\frac{1}{\hbar^{\frac{n}{2}}n!} 
\Gamma\left(\frac{n}{2}+\frac{1}{4}\right)$. 
Now using some properties of the $\Gamma$ function and the formula (2.8) of \cite{diciassette} 
we find an expression denoted by $\mathcal{C}_{1}(n,\hbar)$ multiplied by $\Gamma\left(\frac{1}{4}\right)$.  
It is a fact from the number theory that $\Gamma\left(\frac{1}{4}\right)$ is 
related to the arithmetic--geometric mean denoted by $\agm$, details are in \cite{diciannove} and \cite{venti}.

\hfill $\Box$

\begin{proposition}
Let $\hbar$ be the Plank constant then:

\begin{equation}
\label{lemniscate1}
\int_{0}^{+\infty}e^{-\frac{1}{3\hbar}x^3+\frac{1}{2\hbar}x^2}dx\,=\, 
\mathcal{C}_{2}(n,\hbar)\cdot 
\frac{\pi^{\frac{1}{3}}2^{-\frac{2}{9}}}{3^{-\frac{5}{12}}\agm\left(2,\sqrt{2+\sqrt{3}}\right)^{\frac{1}{3}}}.
\end{equation}
\noindent
where

\begin{equation}
\begin{multlined}[t][12.5cm]
\mathcal{C}_{2}(n,\hbar)\,=\, 
\sum_{n=0}^{+\infty}\frac{1}{(2\hbar)^nn!}3^{\frac{2}{3}n}\hbar^{\frac{2}{3}n}
\prod_{k=1}^{+\infty}\frac{k\left(\frac{2}{3}n+k-1\right)}{\left(k-\frac{1}{3}\right)
  \left(\frac{2}{3}n+k-\frac{2}{3}\right)}\cdot 
  \Gamma\left(\frac{2}{3}n\right).
\end{multlined}
\end{equation}

\end{proposition}

$\Proof.$

\begin{equation}
\begin{multlined}[t][12.5cm]
\int_{0}^{+\infty}e^{-\frac{1}{3\hbar}x^3+\frac{1}{2\hbar}x^2}dx\,=\, \int_{0}^{+\infty}
\left(\sum_{n=0}^{+\infty}\left(\frac{1}{2\hbar}\right)^n\frac{x^{2n}}{n!}\right)e^{-\frac{1}{3\hbar}x^3}dx.
\end{multlined}
\end{equation}

We substitute the variable $\frac{x^3}{3\hbar}\,=\,u$, then we have 

\begin{equation}
\begin{multlined}[t][12.5cm]
\int_{0}^{+\infty}e^{-\frac{1}{3\hbar}x^3+\frac{1}{2\hbar}x^2}dx\\ \,=\, 
\sum_{n=0}^{+\infty}\frac{1}{(2\hbar)^nn!}3^{\frac{2}{3}n-\frac{2}{3}}\hbar^{\frac{2}{3}n+\frac{1}{3}}
\int_{0}^{+\infty}u^{\left(\frac{2}{3}n+\frac{1}{3}\right)-1}e^{-u}du\\ \,=\,
\sum_{n=0}^{+\infty}\frac{1}{(2\hbar)^nn!}3^{\frac{2}{3}n-\frac{2}{3}}\hbar^{\frac{2}{3}n+\frac{1}{3}}
\Gamma\left(\frac{2}{3}n+\frac{1}{3}\right).
\end{multlined}
\end{equation}

Now by the proposition 2.1 of \cite{diciassette} with $x=\frac{2}{3}n$ and $b=\frac{1}{3}$ 
we have an expression for the gamma function:

\begin{equation}
  \begin{multlined}[t][12.5cm]
  \Gamma\left(\frac{2}{3}n+\frac{1}{3}\right)\,=\,
  \prod_{k=1}^{+\infty}\frac{k\left(\frac{2}{3}n+k-1\right)}{\left(k-\frac{1}{3}\right)
  \left(\frac{2}{3}n+k-\frac{2}{3}\right)}\cdot 
  \frac{\Gamma\left(\frac{2}{3}n\right)}{\Gamma\left(\frac{2}{3}\right)}\\
  \,=\,\prod_{k=1}^{+\infty}\frac{k\left(\frac{2}{3}n+k-1\right)}{\left(k-\frac{1}{3}\right)
  \left(\frac{2}{3}n+k-\frac{2}{3}\right)}\cdot 
  \frac{\Gamma\left(\frac{2}{3}n\right)\Gamma\left(\frac{1}{3}\right)3}{2\pi \sqrt{3}}.
  \end{multlined}
\end{equation}

We find an expression denoted by $\mathcal{C}_{2}(n,\hbar)$ multiplied by $\Gamma\left(\frac{1}{3}\right)$.  
As in the previous proposition also the factor $\Gamma\left(\frac{1}{3}\right)$ is 
related to the arithmetic--geometric mean denoted by $\agm$, details are in \cite{diciannove} and 
\cite{venti}.
\hfill $\Box$

\begin{proposition}
Let 
  
$$\psi_{2,4}(v,w)\,=\,-\frac{1}{4\hbar}\|v-w\|^4+\frac{1}{2\hbar}\|v-w\|^2,$$
\noindent
we have that:

\begin{equation}
\label{lemniscate2}
\int_{\mathbb{R}^n}e^{\psi_{2,4}(v,w)}dv\,=\,  
\frac{\pi^{\frac{n}{2}}}{\Gamma\left(\frac{n}{2}+1\right)}\sum_{m=0}\frac{2^{n-4}\hbar^{\frac{n-2m}{4}}}{m!}\Gamma\left(\frac{n+2m}{4}\right).
\end{equation}

\end{proposition}

$\Proof.$
Setting $v-w\,=\,t$ and using the spherical coordinates we have that:

\begin{equation}
\label{lemniscate3}
\begin{multlined}[t][12.5cm]
\int_{\mathbb{R}^n}e^{\psi_{2,4}(v,w)}dv\\
\,=\, \int_{\phi_{n-1}=0}^{2\pi} 
\int_{\phi_{n-2}=0}^{\pi}\cdots  \int_{\phi_{1}=0}^{\pi} \int_{0}^{+\infty}e^{-\frac{1}{4\hbar}r^4+\frac{1}{2\hbar}r^2}
r^{n-1}\\ \cdot\sin^{n-2}{\phi_{1}}\sin^{n-3}{\phi_{2}}\cdots 
\sin{\phi_{n-2}}drd\phi_{1}\cdots d\phi_{n-1}\\
\,=\, 
\int_{S^{n-1}}d_{S^{n-1}}V\int_{0}^{+\infty}r^{n-1}e^{-\frac{1}{4\hbar}r^4+\frac{1}{2\hbar}r^2}dr.
\end{multlined}
\end{equation}

All angular integrals are equal to $V_{n-1}(1)\,=\,\frac{\pi^{\frac{n}{2}}}{\Gamma\left(\frac{n}{2}+1\right)}$ 
and for the second, we use the Taylor expansion of the exponential in this way:

\begin{equation}
  \int_{0}^{+\infty}r^{n-1}\left(\sum_{m=0}^{+\infty}\frac{r^{2m}}{(2\hbar)^mm!}\right) 
  e^{-\frac{1}{4\hbar}r^4}dr,
\end{equation}
\noindent
now substituting $r^4\,=\,4\hbar u$ and, after some calculations, we refind the 
Gamma function:

\begin{equation}
  \int_{0}^{+\infty}r^{n-1}\left(\sum_{m=0}^{+\infty}\frac{r^{2m}}{(2\hbar)^mm!}\right) 
  e^{-\frac{1}{4\hbar}r^4}dr\,=\, 
  \sum_{m=0}\frac{2^{n-4}\hbar^{\frac{n-2m}{4}}}{m!}\Gamma\left(\frac{n+2m}{4}\right).
\end{equation}

\hfill $\Box$

\section{The Berezin--Gaussian integral: a generalization of the Gaussian integral}

In QFT (quantum field theory) Gaussian integrals are treated in the case of anticommutative quantities:

\begin{equation}
\label{anticommutation}
\psi_{1}\psi_{2}\,=\, -\psi_{2}\psi_{1}, \ \ \psi_{1}^2\,=\, \psi_{2}^{2} \,=\, 0,
\end{equation}
\noindent 
called Grassmann numbers or $G$-numbers. We want to treat an integral of the form:

\begin{equation}
\label{antiGaussInt1}
\int_{\Lambda} e^{-\psi_{1}^* r \psi_{1}} d\psi_{1}^*d\psi_{1},
\end{equation}
\noindent 
where $\psi_{1}^*$ is the complex conjugation of a complex Grassmann number with the convention 
that $\left(\psi_{1}\psi_{2}\right)^*=\psi_{2}^*\psi_{1}^*=-\psi_{1}^*\psi_{2}^*$, $r$ is a real exponent 
and $\Lambda$ is the exterior algebra of polynomials in anticommuting elements.

In order to treat integrals of the form $(\ref{antiGaussInt1})$ we observe that the Taylor expansion of 
a function is $f(\psi)= g_{1}+g_{2}\psi$ where $g_{1},g_{2}$ can be both ordinary or Grassmann numbers. 
The integral, called the Berezin integral, is a linear functional with the properties:

\begin{itemize}
 \item[$1)$] $\int_{\Lambda}\frac{\partial f}{\partial \psi}d\psi\,=\, 0$;
 \item[$2)$] $\int_{\Lambda}\psi d\psi  \,=\, 1$;
\end{itemize}
\noindent 
for every $f\in \Lambda$ where $\frac{\partial }{\partial \psi}$ is a left or right derivative. We can treat the $n$ dimensional case $\psi=\left(\psi_{1}, \ldots,\psi_{n}\right)$ in the same manner:

\begin{equation}
\label{antiGaussInt2}
\int_{\Lambda^n} f(\psi) d\psi\,=\, \int_{\Lambda}\left( \cdots \left(\int_{\Lambda} f(\psi) d\psi_{1}\right) \cdots \right)d\psi_{n}.
\end{equation}

The differentiation is defined as:

\begin{equation}
\label{antiGaussInt3}
\frac{\partial}{\partial \psi} \equiv \int_{\Lambda} \cdot \ \ d\psi,
\end{equation}
\noindent 
and this is consistent with previous definitions. 

Now using the properties of the Berezin integral we can treat the integral $(\ref{antiGaussInt1})$ with $\frac{r}{2}$ real:

\begin{equation}
\label{antiGaussInt2}
\int_{\Lambda} e^{-\psi_{1}^* \frac{r}{2} \psi_{1}} d\psi_{1}^*d\psi_{1}\,=\, \int_{\Lambda} \left(1-\psi_{1}^* \frac{r}{2} \psi_{1}\right)d\psi_{1}^*d\psi_{1}\,=\, \frac{r}{2}.
\end{equation}

The result in $n$ dimension with an hermitian matrix $A=\{A_{ij}\}$ is given by:

\begin{equation}
\label{antiGaussInt3}
\begin{multlined}[t][12.5cm]
\left(\int_{\Lambda^n} e^{-\frac{1}{2}\psi_{i}^* A_{ij} \psi_{j}} \prod_{k=1}^{n}d\psi_{k}^*d\psi_{k}\right)^2\\
\,=\,\frac{1}{n!}\int_{\Lambda^n}[-\psi_{i_{1}}^* A_{i_{1}j_{1}} \psi_{j_{1}}]\ldots [-\psi_{i_{n}}^* A_{i_{n}j_{n}} \psi_{j_{n}}]d\psi_{1}^*d\psi_{1}\ldots d\psi_{n}^*d\psi_{n}\\
\,=\,\frac{1}{n!}\int_{\Lambda^n}\left(\psi_{j_{1}}^*\ldots \psi_{j_{n}}^{*}A_{i_{1}j_{1}}\ldots A_{i_{n}j_{n}}\left(\psi_{j_{n}}\ldots \psi_{j_{1}}d\psi_{1}\ldots d\psi_{n}\right)\right)d\psi_{n}^*\ldots d\psi_{1}^*\\ \,=\, \frac{1}{n!}\epsilon_{i_{1}\ldots i_{n}}\epsilon_{j_{1}\ldots j_{n}}A_{i_{1}j_{1}}\cdots A_{i_{n}j_{n}}\,=\, \det{A},
\end{multlined}
\end{equation}
\noindent 
where $\epsilon_{i_{1}\ldots i_{n}}=\int_{\Lambda^{n}}\psi_{i_{1}}\ldots \psi_{i_{n}} d\psi_{n} \ldots d\psi_{1}$ and where there is only one non--zero term in the expansion. So

\begin{equation}
\label{antiGaussInt3}
\begin{multlined}[t][12.5cm]
\int_{\Lambda^n} e^{-\frac{1}{2}\psi_{i}^* A_{ij} \psi_{j}} \prod_{k=1}^{n}d\psi_{k}^*d\psi_{k}\,=\, \sqrt{\det{A}},
\end{multlined}
\end{equation}
\noindent 
that differ from the usual Gaussian integral essentially because $\sqrt{\det{A}}$ is in the numerator instead denominator (details can be founded in \cite{nove}).
This integral is also called the fermionic Gaussian integral.

\subsection{On the Laplace method and a ``mixed Gaussian Integral''}

Now we look if we can say something on integrals of the following form:

\begin{equation}
\label{LantiGaussInt1}
\int_{\Lambda^2} f(\psi,\psi^{*})e^{-\frac{k}{2}g(\psi, \psi^*)} d\psi^*d\psi,
\end{equation}
\noindent 
where we consider the case with $n\,=\,2$. We can apply a Laplace method? 
We can choose $g\,=\,-\psi^* \frac{k}{2} \psi$ an ``heuristic'' approach shows that if $f\,=\, f_{1}+f_{2}\psi+f_{3}\psi^{*}$ and $k\,=\,\frac{1}{\hbar}$ the integral can be computed:

\begin{equation}
\label{LantiGaussInt2}
\int_{\Lambda^2} f(\psi,\psi^{*})e^{-\frac{1}{2\hbar}g(\psi, \psi^*)} d\psi^*d\psi\,=\, \frac{1}{2\hbar}.
\end{equation}

Always with an ``heuristic'' point of view we can consider the following special ``mixed Gaussian integral'':

\begin{equation}
\label{mixedGaussInt1}
\begin{multlined}[t][12.5cm]
\int_{\Lambda^n\times \mathbb{R}^{n}} e^{-\frac{k}{2}\psi_{i}^* \id_{ij} \psi_{j}-\frac{1}{2}k\|x\|^{2}} dx_{1} \cdots dx_{n}\prod_{l=1}^{n}d\psi_{l}^*d\psi_{l},
\end{multlined}
\end{equation}
\noindent 
where $\psi_{i}$ are anticommutative variables and $x_{i}$ are commutative. In this case $A_{ij}\,=\, \id_{ij}$ 
is the identity matrix and $k\,=\,\frac{1}{\hbar}$ is the ``quantum parameter''. 
A combination of tools of previous sections shows that:

\begin{equation}
\label{mixedGaussInt1}
\begin{multlined}[t][12.5cm]
\int_{\Lambda^n\times \mathbb{R}^{n}} e^{-\frac{k}{2}\psi_{i}^* \id_{ij} \psi_{j}-\frac{1}{2}k\|x\|^{2}} dx_{1} \cdots dx_{n}\prod_{l=1}^{n}d\psi_{l}^*d\psi_{l}\,=\, (2\pi)^{\frac{n}{2}},
\end{multlined}
\end{equation}
\noindent
independently from $k$!

\section{Quantum Gaussian integrals with a perturbative term}

In this subsection we consider the presence of a little perturbation term $P(x)$ 
of the original Gaussian integral. We examine the case where the term $P(x)$ is a polynomial 
that consists of a single term:

\begin{equation}
 \label{perturbation}
 P(x)\,=\, \frac{\lambda}{k!}x^k,
\end{equation}
\noindent
where $k\in\mathbb{N}$ and $\lambda >0$. We have the following result in one dimension.

\begin{proposition}
Let $P(x)$ be a perturbation term with the form of $(\ref{perturbation})$, then 

\begin{equation}
  \label{perturbation2}
  \int_{-\infty}^{+\infty}e^{-\frac{1}{2}ax^2+P(x)+Jx}dx\,=\, 
 \left(\frac{2\pi}{a}\right)^{\frac{1}{2}}e^{\frac{\lambda}{k!}\left(\frac{d}{dJ}\right)^k}e^{\frac{J^2}{2a}},
\end{equation}
\noindent 
where $a>0$ and $Jx$ is a source term.
\end{proposition}

$\Proof.$ 

We start with the identity:

\begin{equation}
  \label{perturbation23}
  \int_{-\infty}^{+\infty}e^{-\frac{1}{2}ax^2+Jx}dx\,=\, 
 \left(\frac{2\pi}{a}\right)^{\frac{1}{2}}e^{\frac{J^2}{2a}},
\end{equation}
\noindent
that follows by completing the square and changing coordinates in the integral. 
If we consider in addition the perturbative term, then:

\begin{equation}
  \label{perturbation24}
  \int_{-\infty}^{+\infty}e^{-\frac{1}{2}ax^2+P(x)+Jx}dx\,=\,  
  \int_{-\infty}^{+\infty}\left(\sum_{n=0}^{+\infty}\frac{P(x)^n}{n!}\right)e^{-\frac{1}{2}ax^2+Jx}dx.
\end{equation}

Let us consider this other formula :

\begin{equation}
  \label{perturbation2345}
  \int_{-\infty}^{+\infty}x^ke^{-\frac{1}{2}ax^2+Jx}dx\,=\, 
  \left(\frac{d}{dJ}\right)^k\int_{-\infty}^{+\infty}e^{-\frac{1}{2}ax^2+Jx}dx.  
  \end{equation}

Using $(\ref{perturbation23})$ the equation $(\ref{perturbation2345})$ may be written as:

\begin{equation}
  \label{perturbation2346}
  \int_{-\infty}^{+\infty}x^ke^{-\frac{1}{2}ax^2+Jx}dx\,=\, 
  \left(\frac{2\pi}{a}\right)^{\frac{1}{2}}\left(\frac{d}{dJ}\right)^ke^{\frac{J^2}{2a}}.
 \end{equation}

Now the proposition follows because:

\begin{equation}
  \label{perturbation24}
   \begin{multlined}[t][12.5cm]
     \int_{-\infty}^{+\infty}e^{-\frac{1}{2}ax^2+P(x)+Jx}dx \\ \,=\,  \sum_{n=0}^{+\infty}\frac{1}{n!}\int_{-\infty}^{+\infty}P(x)^ne^{-\frac{1}{2}ax^2+Jx}dx 
   \,=\,\sum_{n=0}^{+\infty}\frac{1}{n!}\left(\frac{\lambda}{k!}\left(\frac{d}{dJ}\right)^k\right)^n\int_{-\infty}^{+\infty}e^{-\frac{1}{2}ax^2+Jx}dx 
   \\
  \,=\,\sum_{n=0}^{+\infty}\frac{1}{n!}\left(\frac{\lambda}{k!}\left(\frac{d}{dJ}\right)^k\right)^n \left(\frac{2\pi}{a}\right)^{\frac{1}{2}}e^{\frac{J^2}{2a}}\,=\, 
  \left(\frac{2\pi}{a}\right)^{\frac{1}{2}}e^{\frac{\lambda}{k!}\left(\frac{d}{dJ}\right)^k}e^{\frac{J^2}{2a}}. 
   \end{multlined}  \end{equation}

\hfill $\Box$

Further information on the perturbative methods in quantum field theory can be found in the book of 
\cite{trentasei}.

\section{Gaussian integrals and the Boys integrals}

In this last section we  will consider a class of Gaussian functions of type: 

$$x^ly^mz^ne^{-ar_{A}^2},$$
\noindent 
where $l,m,n\in \mathbb{Z}_{+}$, $a>0$, $r_{A}^{2}\,=\,x_{A}^2+y_{A}^2+z_{A}^2$, 
$x_{A}\,=\,x-A_{x}$, $y_{A}\,=\, y- A_{y}$, $z_{A}\,=\, z-A_{z}$ and 
$A_{x},A_{y},A_{z}\in\mathbb{R}$. We observe that this class contains the 
complete system of hermitian functions around any point.

\begin{proposition}
  Let $a,b>0$ and $A_{x},A_{y},A_{z},B_{x},B_{y},B_{z}\in\mathbb{R}$, then 
  
\begin{itemize}
  \item[$1)$] 
  
  \begin{equation}\label{boys1}
    \begin{multlined}[t][12.5cm]
    \int_{\mathbb{R}^3}e^{-ar_{A}^2}e^{-br_{B}^2}\,dxdydz \,=\, 
    \left(\frac{\pi}{a+b}\right)^{\frac{3}{2}}e^{-\frac{abS}{a+b}};
  \end{multlined}
  \end{equation}
  
\item[$2)$]

 \begin{equation}\label{boys2}
    \begin{multlined}[t][12.5cm]
    \int_{\mathbb{R}^3}e^{-ar_{A}^2}\left(-\frac{1}{2}\Delta\right)e^{-br_{B}^2}\,dxdydz \,=\, 
    \left(\frac{\pi}{a+b}\right)^{\frac{3}{2}}\left(\frac{3ab}{a+b}-\frac{2Sa^2b^2}{(a+b)^2}\right)e^{-\frac{abS}{a+b}};
  \end{multlined}
  \end{equation}

\item[$3)$]

 \begin{equation}\label{boys3}
    \begin{multlined}[t][12.5cm]
    \int_{\mathbb{R}^3}\frac{e^{-ar_{A}^2}e^{-br_{B}^2}}{r_{C}}\,dxdydz \,=\, 
    \frac{2\pi}{p}F_{0}\left(pR_{CP}^2\right);
  \end{multlined}
  \end{equation}
\end{itemize}\noindent 
where $S\,=\, \sum_{i=x,y,z}(A_{i}-B_{i})^2$, $\Delta$ is the Laplacian in 
$\mathbb{R}^3$, $C$ is the center position of the nucleus, $p\,=\, a+b$, 
$P_{i}\,=\, \frac{aA_{i}+bB_{i}}{p}$ for $i=x,y,z$ called the center of mass, $R_{CP}^2=\sum_{i=x,y,z}(C_{i}-P_{i})^2$ 
and $F_{n}(x)\,=\, \int_{0}^{1}e^{-xt^{2}}t^{2n}\,dt$ is called the ``Boys function''.  

 \end{proposition}

$\Proof.$ The first integral can be written as

 \begin{equation}\label{boys1p}
    \begin{multlined}[t][12.5cm]
    \int_{\mathbb{R}^3}e^{-ar_{A}^2}e^{-br_{B}^2}\,dxdydz \\ \,=\,
    \int_{\mathbb{R}^3}e^{-a(x-A_{x})^2-a(y-A_{y})^2-a(z-A_{z})^2-b(x-B_{x})^2-b(y-B_{y})^2-b(z-B_{z})^2}\,dxdydz,
  \end{multlined}
  \end{equation}
\noindent
we consider only the first of the three integrals and it is equal to:

\begin{equation}
  \label{boys1pp}
  \int_{\mathbb{R}}e^{-(a+b)x^2+2(bB_{x}+aA_{x})x-aA_{x}^{2}-bB_{x}^{2}}\,dx,
\end{equation}
\noindent 
the other two are the same integral but with $y$ and $z$ instead $x$. The 
integral $(\ref{boys1pp})$ can be easily solved completing the square and 
finding that is equal to:

$$ \sqrt{\frac{\pi}{a+b}}e^{-\frac{(A_{x}-B_{x})^2ab}{a+b}},$$
\noindent 
the result for the first integral follows considering the other two integrals.
For the second we reduce the problem in the estimation of:

\begin{equation}\label{boys2p}
    \begin{multlined}[t][12.5cm]
    \int_{\mathbb{R}}e^{-a(x-A_{x})^2}\left(-\frac{1}{2}\frac{\partial^2}{\partial 
    x^2}\right)e^{-b(x-B_{x})^2}\,dx,
  \end{multlined}
  \end{equation}
  \noindent 
  that gives two integrals:
  
  \begin{equation}\label{boys2pp}
    \begin{multlined}[t][12.5cm]
   \,=\,b\int_{\mathbb{R}}e^{-a(x-A_{x})^2-b(x-B_{x})^2}\,dx-2b^2\int_{\mathbb{R}}(x-B_{x})^2e^{-a(x-A_{x})^2-b(x-B_{x})^2}\,dx.
  \end{multlined}
  \end{equation}

The first is the same of $(\ref{boys1pp})$ and the second gives:

$$ \frac{1}{2}\frac{(2A^2+b+2B^2a^2+a-4Ba^2A)e^{-\frac{ab(A_{x}-B_{x})^2}{a+b}}}{(a+b)^2}\cdot 
\sqrt{\frac{\pi}{a+b}}.$$

After algebraic semplifications and recollecting, we have that:

\begin{equation}\label{boys2pp}
    \begin{multlined}[t][12.5cm]
    \int_{\mathbb{R}}e^{-a(x-A_{x})^2}\left(-\frac{1}{2}\frac{\partial^2}{\partial 
    x^2}\right)e^{-b(x-B_{x})^2}\,dx\\
    \,=\, 
    \left[\frac{ab}{a+b}-\frac{2a^2b^2(A_{x}-B_{x})^2}{(a+b)^2}\right]e^{-\frac{ab(A_{x}-B_{x})^2}{a+b}}.
  \end{multlined}
  \end{equation}

We remember the initial reduction and, considering the other two integrals in $y$ 
and $z$, the result follows.

The last integral is a ``Coulomb integral'' that can be solved using a formula for 
the product of two Gaussians. We recall here the formula:

$$ e^{-ax_{A}^2}\cdot e^{-bx_{B}^{2}}\,=\, e^{-\frac{ab}{a+b}X_{AB}^2}\cdot 
e^{-px_{P}^{2}},$$
\noindent 
where $p=a+b$ (is the total exponent), $X_{AB}=(A_{x}-B_{x})$ and $P_{x}\,=\, \frac{aA_{x}+bB_{x}}{p}$. 
The quantity $\frac{ab}{a+b}$ is called the reduced exponent, we will use this 
formula after. Now before to start we consider another identity very easy to 
show using the classical Gaussian integral:

\begin{equation}
 \label{identitycoulomb}
 \frac{1}{r_{C}}\,=\, 
 \frac{1}{\sqrt{\pi}}\int_{-\infty}^{+\infty}e^{-r_{C}^2t^2}dt.
\end{equation}

We transform the original $3$--dimensional integral in the following 
$4$-dimensional integral:

\begin{equation}\label{boys3p}
    \begin{multlined}[t][12.5cm]
    \int_{\mathbb{R}^3}\frac{e^{-ar_{A}^2}e^{-br_{B}^2}}{r_{C}}\,dr \\ 
    \,=\, \frac{1}{\sqrt{\pi}}\int_{\mathbb{R}^4}e^{-pr_{P}^2}e^{-t^2r_{C}^{2}}\,drdt\\
   \,=\, \frac{1}{\sqrt{\pi}}\int_{\mathbb{R}}\int_{\mathbb{R}^3}e^{-(p+t^2)r^2}dr 
   e^{-\frac{pt^2}{p+t^2}R_{CP}^2}\,dt,
  \end{multlined}
  \end{equation}
\noindent 
where the last passage follows from the product formula for Gaussians. Now 
the integral in $dr$ is equal to $\frac{\pi}{p+t^2}$ and the 
initial integral in $3$--dimension is now reduced to be an integral in 
$1$--dimension:

\begin{equation}\label{boys3pp}
    \begin{multlined}[t][12.5cm]
    \int_{\mathbb{R}^3}\frac{e^{-ar_{A}^2}e^{-br_{B}^2}}{r_{C}}\,dr \,=\, 
\frac{2}{\sqrt{\pi}}\int_{0}^{+\infty}\left(\frac{\pi}{p+t^2}\right)e^{-\frac{pt^2}{p+t^2}R_{CP}^2}\,dt.
  \end{multlined}
  \end{equation}

We change variable setting $w^2\,=\,\frac{t^2}{p+t^2}$, so $pt^{-3}dt\,=\, w^{-3}dw$ 
and the integral reduces to 

\begin{equation}\label{boys3ppp}
    \begin{multlined}[t][12.5cm]
    \int_{\mathbb{R}^3}\frac{e^{-ar_{A}^2}e^{-br_{B}^2}}{r_{C}}\,dr \,=\, 
\frac{2\pi}{p}\int_{0}^{1}e^{-pw^2R_{CP}^2}\,dw\,=\frac{2\pi}{p}F_{0}\left(pR_{CP}^2\right).
  \end{multlined}
  \end{equation}

\hfill $\Box$

In \cite{ventidue} the author proves that all Schr\"{o}dinger integrals can be 
estimated for this class of functions, moreover general molecular integrals can be reduced to these integrals. 
A Schr\"{o}dinger integral is a solution of the many--electron Schr\"{o}dinger 
equation for stationary states expressed in atomic units. The atomic units have 
been chosen such that the foundamental electron properties are all equal to one atomic unit. This is the 
reason because the quantum parameter doesn't appears in this section. 
An example of four--center integral is treated in the following letter \cite{ventitre}.

\end{document}